\documentclass[11pt,leqno,twoside]{article}

\usepackage{amsfonts,amsmath,amsthm,amssymb}
\usepackage{enumerate}
\usepackage{graphics,graphicx,subfigure}

\usepackage{color}
\usepackage{todonotes}
\usepackage{cancel}
\usepackage{url}
\usepackage{hyperref}
\usepackage{makeidx}
\usepackage{showidx}
\usepackage{multicol}        % used for the two-column index
\usepackage{xspace}
\usepackage{stmaryrd}        % for \llbracket, \rrbracket
\usepackage{pifont}          % for \ding{227}
\usepackage{fancybox}        % for oval boxes
\usepackage{bm}

 \usepackage{fancyhdr}
\theoremstyle{plain}
 \theoremstyle{definition}
 \newtheorem{lem}{Lemma}
 \newtheorem{defn}[lem]{Definition}
 \newtheorem{thm}[lem]{Theorem}
 \newtheorem{prop}[lem]{Proposition}
 \newtheorem{cor}[lem]{Corollary}
 \newtheorem{notn}[lem]{Notations}
 \newtheorem{pb}[lem]{Problem}
 \newtheorem{form}[lem]{Formulae}
 
 \newtheorem*{rk}{Remark}
 \newtheorem*{com}{Comment}
 \newtheorem*{ex}{Example}
 \theoremstyle{remark}

 \newcommand{\blem}{\begin{lem}}
 \newcommand{\elem}{\end{lem}}
 \newcommand{\bdefn}{\begin{defn}}
 \newcommand{\edefn}{\end{defn}}
 \newcommand{\bthm}{\begin{thm} }
 \newcommand{\ethm}{\end{thm}}
 \newcommand{\bprop}{\begin{prop}}
 \newcommand{\eprop}{\end{prop}}
 \newcommand{\bcor}{\begin{cor}}
 \newcommand{\ecor}{\end{cor}}
 \newcommand{\bnotn}{\begin{notn}}
 \newcommand{\enotn}{\end{notn}}
 \newcommand{\bpb}{\begin{pb}}
 \newcommand{\epb}{\end{pb}}
 \newcommand{\bform}{\begin{form}}
 \newcommand{\eform}{\end{form}}
 \newcommand{\brk}{\begin{rk}}
 \newcommand{\erk}{\end{rk}}
 \newcommand{\bcom}{\begin{com}}
 \newcommand{\ecom}{\end{com}}
 \newcommand{\bex}{\begin{ex}}
 \newcommand{\eex}{\end{ex}}
 \newcommand{\bpf}{\begin{proof}}
 \newcommand{\epf}{\end{proof}}

\newcommand{\cE}{\mathcal{E}}

\newcommand{\cK}{\mathcal{K}}

\newcommand{\cS}{\mathcal{S}}

\newcommand{\bR}{\mathbb{R}}

\newcommand{\be}{\begin{equation}}
\newcommand{\ee}{\end{equation}}
\newcommand{\bal}{\begin{align}}
\newcommand{\eal}{\end{align}}
\newcommand{\ba}{\begin{align*}}
\newcommand{\ea}{\end{align*}}
\newcommand{\bmx}{\begin{matrix}}
\newcommand{\emx}{\end{matrix}}
\newcommand{\bbmx}{\begin{bmatrix}}
\newcommand{\ebmx}{\end{bmatrix}}
\newcommand{\bpmx}{\begin{pmatrix}}
\newcommand{\epmx}{\end{pmatrix}}
\newcommand{\bvmx}{\begin{vmatrix}}
\newcommand{\evmx}{\end{vmatrix}}

\newcommand{\ul}{\underline}

\newcommand{\wh}{\widehat}
\newcommand{\wt}{\widetilde}
\newcommand{\f}{\frac}

\newcommand{\inc}{\subseteq}

\newcommand{\setm}{\setminus}

\newcommand{\tn}[1]{{|\!|\!|}{#1}{|\!|\!|}}
\newcommand{\tnbig}[1]{{\big|\!\big|\!\big|}{#1}{\big|\!\big|\!\big|}}
\newcommand{\tnBig}[1]{{\Big|\!\Big|\!\Big|}{#1}{\Big|\!\Big|\!\Big|}}

\newcommand{\Id}{\mathrm{Id}}

\newcommand{\argmin}{{\rm argmin}\,}

\newcommand{\la}{\lambda}
\newcommand{\La}{\Lambda}
\newcommand{\eps}{\varepsilon}
  
\title{\vspace{-20mm}
Optimal Algorithms for Nonlinear Estimation with Convex Models \rule{166mm}{1.2pt} \vspace{-15mm}}
\author{Simon Foucart\footnote{S. F. is partially supported by a grant from the NSF (DMS-2505204).}  --- Texas A\&M University}
\date{\vspace{-6mm}\rule{100mm}{0.8pt}}

\newcommand\shorttitle{Optimal Algorithms for Nonlinear Estimation with Convex Models}
\newcommand\authors{S. Foucart}

\newcommand{\revX}[1]{{#1}}
\newcommand{\revZ}[1]{{#1}}

\newcommand{\revvY}[1]{{#1}}
\newcommand{\revvZ}[1]{{#1}}
\newcommand{\revvme}[1]{{#1}}

\begin{document}

\maketitle

%% Add abstract, keywords, and AMS classification
\vspace{-15mm}
\begin{abstract}
A linear functional of an object from a convex symmetric set can be optimally estimated,
in a worst-case sense,
by a linear functional of observations made on the object.
This well-known fact is extended here to a nonlinear setting:
other simple functionals of the object can be optimally estimated  by functionals of the observations that share a similar simple structure.
This is established for the maximum of several linear functionals
and even for the $\ell$th largest among them.
Proving the latter requires an unusual \revvme{variation} of the analytical Hahn--Banach theorem.
The existence results are accompanied by practical recipes relying on convex optimization to construct the desired functionals,
thereby justifying the term of estimation algorithms.  
\end{abstract}



\vspace{-5mm}
\begin{center}
\rule{100mm}{0.8pt}
\end{center}

%%%%%%%%%%%%%%%%%
%% The main text starts here %%
%%%%%%%%%%%%%%%%%

\section{Introduction}

Suppose that an unknown object $f$,
living in a vector space $F$,
is acquired via a vector $\La(f) \in \bR^m$ of linear observations $\la_1(f),\ldots,\la_m(f) \in \bR$,
in addition to the prior knowledge that it belongs to a prescribed model set $\cK$.
Estimating a quantity \revvZ{$\gamma(f)$},
living in some normed space $X$,
boils down to devising a recovery map \revvZ{$\delta: \bR^m \to X$}.
From a worst-case perspective,
its performance is measured through
\revvZ{$$
e(\delta) := \sup_{f \in \cK} \| \gamma(f) -\Delta(\La(f)) \|.
$$
}For a fixed \revZ{observation map} $\La$,
constructing a map \revvZ{$\delta$} that minimizes the above
is the essence of Optimal Recovery~\cite{MicRiv}, 
a lead-in to the field of Information-Based Complexity \cite{TWW}, 
where \revZ{one is also interested in minimizing over $\La$,
with possibly some restrictions on the nature of the observations
(see the comprehensive survey \cite{KriUll} for results in this direction).
From here on,
the Optimal Recovery framework prevails, so the dependence on the linear observation map $\La: F \to \bR^m$---considered a fixed entity---will not be made explicit.
For instance,
it is implicit in the lower bound $e_\flat$ below,
asserting that
for any \revvZ{$\delta: \bR^m \to X$},}
\revvZ{$$
e(\delta) \ge e_\flat,
\qquad \mbox{where} \quad
e_\flat := \sup_{\substack{f',f'' \in \cK \\ \La(f') = \La(f'')}} \f{\| \gamma(f') - \gamma(f'') \|}
{2}.
$$
}\revZ{This} is well known and easy to derive,
and so is the complementary inequality \revvZ{$e(\gamma \circ \Delta^{\rm cons}) \le 2 e_\flat$},
where $\Delta^{\rm cons}: \bR^m \to F$ is any
data- and model-consistent map\footnote{\revX{A data- and model-consistent map is usually called an interpolatory algorithm in Information-Based Complexity,
see e.g. \cite[p,49--53]{TWW}.}},
i.e.,
one satisfying $\La(\Delta^{\rm cons}(y))= y$ and $\Delta^{\rm cons}(y) \in \cK$ for all $y \in \La(\cK)$.
In short, the map \revvZ{$\gamma \circ \Delta^{\rm cons}$}
is nearly optimal (with a factor $2$).
To surpass this, one seeks maps \revvZ{$\delta^{\rm opti}: \bR^m \to X$} yielding \revvZ{$e(\delta^{\rm opti}) \le e_\flat$}, which are therefore genuinely optimal (with a factor $1$).
This is achievable \revvZ{when $\gamma$ is a linear map taking values in $X = \bR$,
i.e., it is a linear functional.}
%when $\Gamma = \gamma$ is a linear functional---throughout, lowercase Greek letters are used for functionals. 
\revvZ{Such a} classical result is due to Smolyak (see e.g.  \cite[Theorem~4.7]{NovWoz}) when the set $\cK$ is convex and symmetric about the origin
and was refined by Sukharev \cite{Suk} when $\cK$ is merely convex.
These two results also guarantee the existence of a genuinely optimal estimation \revvme{functional} $\revZ{\delta^{\rm opti}}: \bR^m \to \bR$ which is linear in the former case and affine in the latter \revZ{case}.

The goal of this article is, firstly, to exhibit examples of nonlinear yet `simple' functionals $\gamma: F \to \bR$
for which 
\revZ{there exist `simple' estimation \revvme{functionals} $\delta^{\rm opti}: \bR^m \to \bR$ that are genuinely optimal, in the sense that
$$
e(\delta^{\rm opti})
\le e(\delta)
\qquad \mbox{for any } \delta: \bR^m \to \bR,
$$
and} secondly, to complement existence results with practical constructions whenever possible. 
As a prototypical example of nonlinear functionals,
the maximum of a function $f$ was considered early on, \revX{often with near optimality in mind (see e.g. \cite{Was}), but not always (see e.g. \cite[p.28--31]{Nov})}.
With genuine optimality in mind,
it has \revX{also} been treated \revX{recently} relatively to a model set $\cK$ consisting of H\"older functions within the space $F$ of continuous functions, see~\cite{FouMax}.
The arguments were specific to the case at hand,
so the situation where the function space $F$ is a reproducing kernel Hilbert space, say,
was not covered.
The results obtained in the present article do apply to such a situation:
for instance, Theorem~\ref{Thm4Max} will guarantee the existence of a genuinely optimal estimation functional $\delta^{\rm opti}: \bR^m \to \bR$ which is convex
and Proposition~\ref{ThmConstruction} will later provide a practical recipe to construct it,
under the proviso that the maximum is over a finite set in order for a convex optimization program to be solvable.
As an extension,
Theorem~\ref{ThmExt} will establish e.g. that the functional outputting the $\ell$th largest value among a number of linear functionals
comes with a genuinely optimal estimation functional $\delta^{\rm opti}: \bR^m \to \bR$ which is the supremum of infima of affine functionals.
This result is based an a \revvme{variation} of the Hahn--Banach dominated extension theorem, see Lemma~\ref{LemHB}.

Here is a brief outline of the organization of the article:
theoretical results such as Theorem~\ref{Thm4Max} and Theorem~\ref{ThmExt} are proved right below (Section~\ref{SecTheory}),
before computational considerations are addressed (Section~\ref{SecComp}),
first at an abstract level and then specialized to two particular examples.

\section{Optimal Estimation of Some Nonlinear Functionals}
\label{SecTheory}

This section starts by uncovering the `simple' form of an optimal functional for the estimation of a nonlinear functional $\gamma$ which appears as the supremum of linear functionals.
It continues by enlarging the result to include $\gamma$'s that are a mixed supremum-infimum of linear functionals.
In the process,
%a seemingly novel refinement 
\revvme{an unusual variation}
of Hahn--Banach dominated extension theorem is
\revvme{exploited.}
% established.
The section ends with a discussion on the handling of possible errors in the observation procedure.

\subsection{The supremum of linear functionals}

In the result presented below,
valid in an arbitrary normed space $F$,
the nonlinear functional $\gamma$ to be optimally estimated is the supremum of linear functionals indexed by a set $I$.
Since the latter need not be finite,
it transpires  that any norm $\revX{ \|f \|} = \max\{ \eta(f): \eta \in F^* \mbox{ with } \revX{\| \eta \|_*} =1 \}$ can be optimally estimated by way of a functional $\delta^{\rm opti}: \bR^m \to \bR$ which,
as the supremum of affine functionals, is a convex functional.
Of course, an analog result could be obtained if $\gamma$ was the infimum of linear functionals
by simply applying Theorem~\ref{Thm4Max} to $-\gamma$.

\bthm
\label{Thm4Max}
Suppose that the model set $\cK$ is convex and contains the origin and that the functional $\gamma: F \to \bR$ is sup-linear,
i.e., of the form
$$
\gamma(f) = \sup_{i \in I} \gamma_i(f),
\qquad f \in F,
$$
where the $\gamma_i \in F^*$ are linear functionals.
Then there exists an optimal estimation functional $\delta^{\rm opti} : \bR^m \to \bR$ which is also sup-affine,
i.e., of the form
$$
\delta^{\rm opti}(y) = \sup_{i \in I} 
\bigg( c^{(i)}_0 + \sum_{k = 1}^m c^{(i)}_k y_k  \bigg),
\qquad y \in \bR^m.
$$
\ethm

\bpf
As a first step, the previously mentioned lower bound is rewritten as
$$
e_\flat = \sup_{\substack{f',f'' \in \cK \\ \La(f') = \La(f'')}} \f{\gamma(f')-\gamma(f'')}{2}
= \sup_{\bar{f} \in \bar{\cK} \cap \ker(\bar{\La})}
\f{\bar{\gamma}'(\bar{f}) - \bar{\gamma}''(\bar{f}) }{2},
$$
where one has introduced the set $\bar{\cK} := \cK \times \cK \inc \bar{F} := F \times F$,
as well as the linear map $\bar{\La}: \bar{F} \to \bR^m$ and the (nonlinear) functionals $\bar{\gamma}',\bar{\gamma}'' : \bar{F} \to \bR$ defined via
$$
\bar{\La}([f';f'']) = \La(f') - \La(f''),
\qquad \bar{\gamma}'([f';f'']) = \gamma(f'),
\qquad \bar{\gamma}''([f';f'']) = \gamma(f'').
$$
With $\rho_{\bar{\cK}}$ denoting the Minkowski functional of $\bar{\cK}$---which is a sublinear (i.e., positively homogeneous and convex) functional 
by virtue of the fact that $\bar{\cK}$ is convex and contains the origin---one has  
$\bar{\gamma}'(\bar{f}) - \bar{\gamma}''(\bar{f}) \le 2 e_\flat \rho_{\bar{\cK}}(\bar{f})$ whenever $\bar{f} \in \ker(\bar{\La})$.
Also introducing linear functionals $\bar{\gamma}'_i, \bar{\gamma}''_i: \bar{F} \to \bR$ defined for all $i \in I$ via
$\bar{\gamma}'_i([f';f'']) = \gamma_i(f')$
and $\bar{\gamma}''_i([f';f'']) = \gamma_i(f'')$,
so that $\bar{\gamma}'(\bar{f}) = \sup_{i \in I} \bar{\gamma}'_i(\bar{f})$,
the previous inequality ensures that, for all $i_\star \in I$,
one has
$\bar{\gamma}'_{i_\star}(\bar{f}) \le 2 e_\flat \rho_{\bar{\cK}}(\bar{f}) + \bar{\gamma}''(\bar{f}) $
whenever $\bar{f} \in \ker(\bar{\La})$.
Fixing $i_\star$ and noticing that $2 e_\flat \rho_{\bar{\cK}} + \bar{\gamma}''$ is a sublinear functional,
the Hahn--Banach dominated extension theorem
(see e.g. \revvme{\cite[Theorem 3.2]{Rud}})
%\cite[Theorem 4, p. 49]{Bol})
applies to guarantee the existence of a linear functional $\mu_{i_\star} \in \bar{F}^*$ such that 
${\mu_{i_{\star}}}_{| \ker(\bar{\La})} = {\bar{\gamma}'_{i_{\star}}}{}_{| \ker(\bar{\La})}$
and $\mu_{i_\star}(\bar{f}) \le 2 e_\flat \rho_{\bar{\cK}}(\bar{f}) + \bar{\gamma}''(\bar{f}) $  for all $\bar{f} \in \bar{F}$.
In view of the equivalence between the fact that $\bar{\gamma}'_{i_{\star}} - \mu_{i_{\star}} $ vanishes on $\ker(\bar{\La})$ and the existence of $c^{(i_\star)} \in \bR^m$ such that
$\bar{\gamma}'_{i_{\star}} - \mu_{i_{\star}} = \sum_{k=1}^m c^{(i_\star)}_k \bar{\la}_k$
\revvme{(see e.g. \cite[Lemma~3.9]{Rud})}, one derives that
$$
\bar{\gamma}'_{i_{\star}}(\bar{f})
- \sum_{k = 1}^m c^{(i_\star)}_k \bar{\la}_k(\bar{f})
\le 2 e_\flat \rho_{\bar{\cK}}(\bar{f}) + \bar{\gamma}''(\bar{f}) 
\qquad \mbox{for all }\bar{f} \in \bar{F}.
$$
In other words,
for any $i_\star \in I$, it holds that
$$
\gamma_{i_\star}(f') 
- \sum_{k = 1}^m c^{(i_\star)}_k
\big( \la_k(f') - \la_k(f'') \big)
\le 2 e_\flat \max\{ \rho_{\cK}(f'), \rho_{\cK}(f'')  \}+ \gamma(f'') 
\qquad \mbox{for all }f',f'' \in F,
$$
which immediately yields the key observation that
\be
\label{KEY}
\bigg( \gamma_{i_\star}(f') 
- \sum_{k = 1}^m c^{(i_\star)}_k \la_k(f') \bigg)
- \bigg(
\gamma(f'') - \sum_{k = 1}^m c^{(i_\star)}_k \la_k(f'')
\bigg)
\le 2 e_\flat 
\qquad \mbox{for all }f',f'' \in \cK.
\ee
At this point, for each $i_\star \in I$, one defines
\be
\label{DefC0}
c^{(i_\star)}_0 = \inf_{f \in \cK} \bigg( 
\gamma(f) - \sum_{k=1}^m c^{(i_\star)}_k \la_k (f) + e_\flat
 \bigg).
\ee
This definition readily implies that,
for each $i_\star \in I$,
one has
$\gamma(f) - \Big( c^{(i_\star)}_0 +  \sum_{k=1}^m c^{(i_\star)}_k \la_k (f) \Big) \ge  - e_\flat$ for all $f \in \cK$.
One can now take the supremum over $i_\star \in I$ to arrive at
\be
\label{KEY1}
\gamma(f) - \sup_{i \in I} \bigg( c^{(i)}_0 +  \sum_{k=1}^m c^{(i)}_k \la_k (f) \bigg) \ge  - e_\flat
\qquad \mbox{for all } f \in \cK.
\ee
Next,
for each $i_\star \in I$,
one selects $f'' \in \cK$ such that $
\gamma(f'') - \sum_{k=1}^m c^{(i_\star)}_k \la_k (f'') + e_\flat$ is equal (or is arbitrarily close) to $c^{(i_\star)}_0$.
Writing \eqref{KEY} for this $f''$ and for $f'$ being an arbitrary $f \in \cK$ leads to
$\gamma_{i_\star}(f) \le e_\flat
+ \Big( c^{(i_\star)}_0 + \sum_{k = 1}^m c^{(i_\star)}_k \la_k(f) \Big)$,
and hence
$\gamma_{i_\star}(f) \le e_\flat
+ \sup_{i \in I}\Big( c^{(i)}_0 + \sum_{k = 1}^m c^{(i)}_k \la_k(f) \Big)$.
Taking the supremum over $i_\star \in I$ yields
\be
\label{KEY2}
\gamma(f) - \sup_{i \in I} \bigg( c^{(i)}_0 +  \sum_{k=1}^m c^{(i)}_k \la_k (f) \bigg) \le   e_\flat
\qquad \mbox{for all } f \in \cK.
\ee
Now defining the estimation functional $\delta^{\rm opti}: \bR^m \to \bR$ by $\delta^{\rm opti}(y) = \sup_{i \in I} \Big( c^{(i)}_0 +  \sum_{k=1}^m c^{(i)}_k y_k \Big)$
\revvme{and}
combining inequalities \eqref{KEY1} and \eqref{KEY2}  show that
$$
e(\delta^{\rm opti}) = \sup_{f \in \cK}
\big|
\gamma(f) - \delta^{\rm opti}(\La(f))
\big| 
\le e_\flat.
$$
This proves the genuine optimality of the sup-affine functional $\delta^{\rm opti}$.
\epf

\revvY{
\brk
Although natural from an Optimal Recovery point of view,
invoking Hahn--Banach dominated extension theorem may appear quite abstract,
so building intuition via an alternative route could be welcome.
Here is a sketch of such a convex-duality route,
pointed out by a referee.
Start by writing $2 e_\flat = \sup\{ \gamma(f') - \gamma(f''): \, f',f'' \in \cK, \, \La(f') = \La(f'') \}$ as
$$
2 e_\flat = \sup_{i \in I} \mu_i,
\qquad 
\mbox{where}
\quad
\mu_i:=
\sup_{\substack{f',f'' \in \cK \\ \La(f') = \La(f'')}} \big( \gamma_i(f')-\gamma(f'') \big).
$$
Since the program defining $-\mu_i$ is convex, consider its dual program 
to obtain (after some work, in particular to justify strong duality)
$$
\mu_i=
\inf_{c^{(i)} \in \bR^m}
\sup_{f',f'' \in \cK }
\bigg(
 \gamma_i(f')-\gamma(f'')
- \sum_{k=1}^m c^{(i)}_k (\la_k(f') - \la_k(f''))
\bigg).
$$
With $c^{(i)}$ representing for simplicity a minimizer of the above,
observe that
$$
\mu_i - e_\flat =
\sup_{f' \in \cK} \bigg( \gamma_i(f') - \sum_{k=1}^m c^{(i)}_k \la_k(f') \bigg)
+ \sup_{f'' \in \cK} \bigg( -\gamma(f'') + \sum_{k=1}^m c^{(i)}_k \la_k(f'') - e_\flat \bigg).
$$ 
Denote the latter supremum as $-c^{(i)}_0$,
in accordance with~\eqref{DefC0}, and deduce that
$$
e_\flat = \sup_{i \in I} (\mu_i - e_\flat) 
= \sup_{f' \in \cK} \sup_{i \in I} \bigg(
\gamma_i(f')  -  \sum_{k=1}^m c^{(i)}_k \la_k(f') - c^{(i)}_0 \bigg)
\ge \sup_{f' \in \cK} \Big( \gamma(f') - \delta^{\rm opti}(\La f') \Big)
$$
where $\delta^{\rm opti}(y) = \sup_{i \in I} \big( \sum_{k=1}^m c^{(i)}_k y_k + c^{(i)}_0 \big)$.
This inequality is one of two that explain the optimality of the estimation functional $\delta^{\rm opti}$.
The other one, namely $\sup_{f'' \in \cK} \big(
- \gamma(f'') + \delta^{\rm opti}(\La f'')   \big) \le e_\flat$,
is a direct consequence of the definition of $c_0^{(i)}$.
\erk
}

\subsection{The mixed supremum-infimum of linear functionals}

It is now acquired that there is a `simple' optimal algorithm for the estimation of the (pointwise) largest value of linear functionals $\gamma_1,\ldots,\gamma_d$, say.
Is there also a `simple' optimal algorithm for the estimation of the (pointwise) $\ell$th largest value of these linear functionals?
The answer is yes, as established in Theorem~\ref{ThmExt} below.
`Simple' in this case means the mixed sup-inf of affine functionals.
This mirrors the sup-inf/inf-sup nature of the functional $\gamma^{\downarrow}_\ell$ associating to $f \in F$
the $\ell$th largest value of $\gamma_1(f),\ldots,\gamma_d(f)$,
namely
$$
\gamma^{\downarrow}_\ell(f)
= \sup_{|I|=\ell} \inf_{i \in I} \gamma_i(f)
= \inf_{|J|=d+1-\ell} \sup_{j \in J} \gamma_j(f).
$$  
The essential feature is to have both a sup-inf and an inf-sup representation of the nonlinear functional at hand, 
as e.g. in Courant--Fischer characterization of the $\ell$th largest eigenvalue. 
Thus, the result will also apply to the difference of suprema of linear functionals,
since
$$
\sup_{i \in I} \mu_i(f) - \sup_{j \in J} \nu_j(f)
= \sup_{i \in I} \inf_{j \in J} (\mu_i - \nu_j)(f)
= \inf_{j \in J} \sup_{i \in I} (\mu_i - \nu_j)(f).
$$
\revvY{A crucial part of the argument will
rely on the following unusual, yet simple, variation on Hahn--Banach dominated extension theorem.}

\blem
\label{LemHB}
Let $V$ be a vector space 
and let $U$ be a finite codimensional subspace of $V$ given as $U = \ker(\eta_1) \cap \cdots \cap \ker(\eta_m)$ for some linear functionals $\eta_1,\ldots,\eta_m \in V^*$.
If $\mu_i \in V^*$, $i \in I$,
are linear functionals on $V$ and if $\rho$ is a sublinear functional on $V$ such that
$$
\inf_{i \in I} \mu_i(u) \le \rho(u)
\qquad \mbox{for all } u \in U,
$$
then there exist scalars $c_1,\ldots,c_m \in \bR$ such that
$$
\inf_{i \in I} \mu_i(v) +
\sum_{k=1}^m c_k \eta_k(v)
\le \rho(v)
\qquad \mbox{for all } v \in V.
$$
\elem

\revvY{
\bpf
The assumption says that the zero linear functional is dominated on the subspace $U$ by the sublinear functional $\rho - \inf_{i \in I} \mu_i$.
Applying the standard Hahn--Banach dominated extension theorem then yields a linear functional $\omega$ such that 
$\omega \le \rho - \inf_{i \in I} \mu_i$ on the whole of $V$ and $\omega_{|U} = 0$.
Since the latter is equivalent
to the existence of scalars $c_1,\ldots,c_m \in \bR$ such that   
$\omega = \sum_{k=1}^m c_k \eta_k$ 
(see e.g. \cite[Lemma~3.9]{Rud}),
the announced result follows. 
\epf
}

With this technical lemma at hand,
the main result of the subsection can now be stated and proved.

\bthm
\label{ThmExt}
Suppose that the model set $\cK$ is convex and contains the origin and that the functional $\gamma: F \to \bR$ can be expressed as
$$
\gamma(f) = \sup_{a \in A} \inf_{i \in I_a} \gamma_i(f)
= \inf_{b \in B} \sup_{j \in J_b} \gamma_j(f),
\qquad f \in F,
$$
where the $\gamma_i \in F^*$ are linear functionals.
Then there exists an optimal estimation functional $\delta^{\rm opti} : \bR^m \to \bR$ which has the form
$$
\delta^{\rm opti}(y) = \sup_{a \in A}
\inf_{b \in B} 
\bigg( c^{(a,b)}_0 + \sum_{k = 1}^m c^{(a,b)}_k y_k  \bigg),
\qquad y \in \bR^m.
$$
\ethm

\bpf
One starts from the lower bound $e_\flat$ expressed once again as
$$
e_\flat = \sup_{\substack{f',f'' \in \cK \\ \La(f') = \La(f'')}} \f{\gamma(f')-\gamma(f'')}{2}
= \sup_{\bar{f} \in \bar{\cK} \cap \ker(\bar{\La})}
\f{\bar{\gamma}'(\bar{f}) - \bar{\gamma}''(\bar{f}) }{2},
$$
where it is recalled that the `bar'-notation  involves
$\bar{\cK} := \cK \times \cK \inc \bar{F} := F \times F$, which has Minkowski functional $\rho_{\bar{\cK}}$,
and  $\bar{\La}: \bar{F} \to \bR^m$, $\bar{\gamma}',\bar{\gamma}'',\bar{\gamma}'_i, \bar{\gamma}''_i : \bar{F} \to \bR$, which are defined via
\begin{align*}
\bar{\La}([f';f'']) = \La(f') - \La(f''),
& & \bar{\gamma}'([f';f'']) = \gamma(f'),
\quad \qquad \bar{\gamma}''([f';f'']) = \gamma(f''), \;\\
& & \bar{\gamma}'_i([f';f'']) = \gamma_i(f'),
\quad \quad \; \; \; \bar{\gamma}''_i([f';f'']) = \gamma_i(f'').
\end{align*}
In this way,  one has 
$\bar{\gamma}'(\bar{f}) 
\le 2 e_\flat \rho_{\bar{\cK}}(\bar{f}) +
\bar{\gamma}''(\bar{f})$
whenever
$\bar{f} \in \ker(\bar{\La})$.
Therefore, for all $a_\star \in A$ 
and all $b_\star \in B$,
$$
\inf_{i \in I_{a_\star}} \bar{\gamma}_i'(\bar{f}) \le 2 e_\flat \rho_{\bar{\cK}}(\bar{f}) + 
\sup_{j \in J_{b_\star}} \bar{\gamma}_j''(\bar{f})
\qquad \mbox{whenever }
\bar{f} \in \ker(\bar{\La}).
$$
This is the time to invoke Lemma \ref{LemHB},  justifying the existence of $c^{(a_\star,b_\star)} \in \bR^m$ such that
$$
\inf_{i \in I_{a_\star}} \bar{\gamma}_i'(\bar{f}) 
-\sum_{k=1}^m c^{(a_\star,b_\star)}_k \bar{\la}_k(\bar{f}) 
\le 2 e_\flat \rho_{\bar{\cK}}(\bar{f}) + 
\sup_{j \in J_{b_\star}} \bar{\gamma}_j''(\bar{f})
\qquad \mbox{for all }
\bar{f} \in \bar{F}.
$$
It follows that,
for any $a_\star \in A$ and $b_\star \in B$, 
one has
\be
\label{KEYInfSup}
\bigg( \inf_{i \in I_{a_\star}} \gamma_i(f') 
-\sum_{k=1}^m c^{(a_\star,b_\star)}_k \la_k(f') \bigg)
- \bigg( \sup_{j \in J_{b_\star}} \gamma_j(f'')
-\sum_{k=1}^m c^{(a_\star,b_\star)}_k \la_k(f'') \bigg)
\le 2 e_\flat  
\; \; \, \mbox{for all }
f',f'' \in \cK.
\ee
At this point, for each $a_\star \in A$ and $b_\star \in B$, one defines
$$
c^{(a_\star,b_\star)}_0
:= 
\inf_{f \in \cK}
\bigg( \sup_{j \in J_{b_\star}} \gamma_j(f)
-\sum_{k=1}^m c^{(a_\star,b_\star)}_k \la_k(f) + e_\flat \bigg).
$$
For each $a_\star \in A$ and $b_\star \in B$,
one immediately obtains that,
for all $f \in \cK$,
$$
\sup_{j \in J_{b_\star}} \gamma_j(f) \ge 
-e_\flat
+ \bigg( c^{(a_\star,b_\star)}_0 + \sum_{k=1}^m c^{(a_\star,b_\star)}_k \la_k(f)  \bigg)
\ge 
-e_\flat
+ \inf_{b \in B} \bigg( c^{(a_\star,b)}_0 + \sum_{k=1}^m c^{(a_\star,b)}_k \la_k(f)  \bigg).
$$
Now, taking the infimum over $b_\star \in B$ on the left-hand side yields
$$
\gamma(f) \ge 
-e_\flat
+ \inf_{b \in B} \bigg( c^{(a_\star,b)}_0 + \sum_{k=1}^m c^{(a_\star,b)}_k \la_k(f)  \bigg)
\qquad \mbox{for all } f \in \cK,
$$
and finally taking the supremum over $a_\star \in A$ on the right-hand side and rearranging gives
\be
\label{KEYInfSup1}
\gamma(f) - \sup_{a \in A} \inf_{b \in B} \bigg( c^{(a,b)}_0 + \sum_{k=1}^m c^{(a,b)}_k \la_k(f)  \bigg)
\ge 
-e_\flat
\qquad \mbox{for all } f \in \cK.
\ee
Next,
for each $a_\star \in A$ and $b_\star \in B$,
one selects $f'' \in \cK$ 
achieving (or coming arbitrarily close to) the infimum defining $c_0^{(a_\star,b_\star)}$,
so writing \eqref{KEYInfSup}
for this $f''$ and for $f'$ being an arbitrary $f \in \cK$ leads to
$$
\inf_{i \in I_{a_\star}} \gamma_i(f)
- \bigg( c^{(a_\star,b_\star)}_0 + \sum_{k=1}^m c^{(a_\star,b_\star)}_k \la_k(f)  \bigg)
\le e_\flat,
$$
from where it follows,
by taking the infimum over $b_\star \in B$, that
$$
\inf_{i \in I_{a_\star}} \gamma_i(f) 
\le e_\flat + 
\inf_{b \in B} \bigg( c^{(a_\star,b)}_0 + \sum_{k=1}^m c^{(a_\star,b)}_k \la_k(f)  \bigg)
\le e_\flat + 
\sup_{a \in A} \inf_{b \in B} \bigg( c^{(a,b)}_0 + \sum_{k=1}^m c^{(a,b)}_k \la_k(f)  \bigg).
$$
Finally, taking the supremum over $a_\star \in A$ leads to
\be
\label{KEYInfSup2}
\gamma(f) - \sup_{a \in A} \inf_{b \in B} \bigg( c^{(a,b)}_0 + \sum_{k=1}^m c^{(a,b)}_k \la_k(f)  \bigg)
\le e_\flat
\qquad \mbox{for all } f \in \cK.
\ee
Defining the estimation functional 
$\delta^{\rm opti}: \bR^m \to \bR$ by $\delta^{\rm opti}(y) = \sup_{a \in A}
\inf_{b \in B} 
\big( c^{(a,b)}_0 + \sum_{k = 1}^m c^{(a,b)}_k y_k  \big)$
and combining inequalities \eqref{KEYInfSup1} and \eqref{KEYInfSup2} show that
$$
e(\delta^{\rm opti}) = \sup_{f \in \cK}
\big|
\gamma(f) - \delta^{\rm opti}(\La(f))
\big| 
\le e_\flat.
$$
This proves the optimality of the sup-inf-affine functional $\delta^{\rm opti}$.
\epf

\subsection{Handling observation error}

In realistic situations,
the observations made on the unknown object $f$ are not exact.
They typically take the form $\la_k(f) + e_k$ for some unknown $e \in \bR^m$.
This vector can be modeled stochastically,
as done in \cite{Don} but not considered here,
or deterministically via an assumption $e \in \cE$, as considered here, where $\cE$ is convex set  containing the origin.
 It is folklore to remark that this `inaccurate scenario' \revvY{formally}
 reduces to the `accurate scenario'
 by focusing on the compound object $(f,e)$
 belonging to the model set $\cK \times \cE$
 and acquired via the linear observations $\wt{\la}_k((f,e)) = \la_k(f) + e_k$.
For instance,
since a sup-linear functional $\gamma$ acting on $f$ is also a sup-linear functional acting on $(f,e)$,
Theorem~\ref{Thm4Max} will still guarantee that there exists \revvme{a sup-affine functional which is optimal for the estimation of $\gamma$.}
\revvY{Here is a more explicit statement for the counterpart of Theorem~\ref{ThmExt}.
\bprop
\label{PropInac}
Suppose that the sets $\cK \inc F$ and $\cE \inc \bR^m$ are convex and contains the origins
and that the functional $\gamma: F \to \bR$ can be expressed as
$$
\gamma(f) = \sup_{a \in A} \inf_{i \in I_a} \gamma_i(f)
= \inf_{b \in B} \sup_{j \in J_b} \gamma_j(f),
\qquad f \in F,
$$
where the $\gamma_i \in F^*$ are linear functionals.
Then there exists a functional $\delta^{\rm opti} : \bR^m \to \bR$ of the form
$$
\delta^{\rm opti}(y) = \sup_{a \in A}
\inf_{b \in B} 
\bigg( c^{(a,b)}_0 + \sum_{k = 1}^m c^{(a,b)}_k y_k  \bigg),
\qquad y \in \bR^m,
$$
which is optimal for the worst-case estimation of $\gamma$, in the sense that, for all $\delta: \bR^m \to \bR$,
$$
\sup_{\substack{f \in \cK \\ e \in \cE}}
|\gamma(f) - \delta^{\rm opti}(\La(f) + e)|
\le \sup_{\substack{f \in \cK \\ e \in \cE}}
|\gamma(f) - \delta(\La(f) + e)|.
$$
\eprop

The practical construction of the optimal estimation functional is addressed in the next section---thanks to the formal reduction, it will be enough to consider the `accurate scenario' there.
A central role shall be played by
the support function of the set $\cK$,
as defined by
$$
\tn{\eta}_\cK := \sup_{f \in \cK} \eta(f),
\qquad \eta \in F^*.
$$ 
In the `inaccurate scenario',
one should use the support function of the set $\cK \times \cE$ instead.
To fix the ideas, suppose as an example that $\cE = \eps B_p^m$ is an origin-centered ball of radius $\eps$ in $\ell_p$-space for some $p \in [1,\infty]$ with conjugate exponent $p' \in [1,\infty]$.
Then, the observation that
$$
\tnBig{ \sum\nolimits_{k=1}^m c_k \wt{\la}_k }_{\cK \times \eps B_p^m}
= \sup_{(f,e) \in \cK \times \eps B_p^m}  \sum\nolimits_{k=1}^m c_k \big( \la_k(f) + e_k \big)
= \tnBig{ \sum\nolimits_{k=1}^m c_k \la_k }_{\cK}
+ \eps \|c\|_{p'}
$$
enables the application of the forthcoming Proposition~\ref{ThmConstruction2}, say,
to determine the coefficients $c^{(a,b)}_k$ of the optimal estimation functional $\delta^{\rm opti}$.

While an analysis of the dependence of these coefficients on the level $\eps$ of observation error is beyond the scope of this article,
a rough estimate of the dependence on $\eps$ of the minimal worst-case error is not hard to come by.
Precisely,
with 
$e_\eps(\delta) := \sup\{ |\gamma(f) - \delta(\La f + e)|: \, f \in \cK, \, e \in \eps B_p^m \}$
denoting the worst-case error of an arbitrary $\delta: \bR^m \to \bR$ at level $\eps$,
one has
$$
\inf_\delta e_\eps(\delta)
\le 
\inf_\delta
e_{\zeta}(\delta)
\le 
\inf_\delta e_\eps(\delta)
+ C_\eps \, (\zeta - \eps)
\qquad \mbox{for any } \zeta \ge \eps \ge 0,
$$
where the constant $C_\eps$ can be taken as the largest $\ell_{p'}$-norm of the coefficient vectors $c^{(a,b)} \in \bR^m$ appearing in the expression of the optimal estimation functional $\delta^{\rm opti}$ at level $\eps$.
Indeed, 
the leftmost inequality follows from the simple fact that $e_\eps(\delta) \le e_{\zeta}(\delta)$ for any $\delta: \bR^m \to \bR$.
As for the rightmost inequality,
it relies on
a fact easily verifiable from the Lipschitz properties of $\inf$ and $\sup$,
namely from $|\delta^{\rm opti}(y) - \delta^{\rm opti}(y')| \le \sup_{a \in A} \sup_{b\in B} \big| \sum_{k=1}^m c^{(a,b)}_k (y_k - y'_k) \Big|$,
as used in 
\begin{align*}
\inf_\delta e_{\zeta}(\delta)
& \le e_{\zeta}(\delta^{\rm opti})
= \sup_{\substack{f \in \cK \\ \|e\|_p \le \zeta }}
|\gamma(f) - \delta^{\rm opti}(\La(f)+e)|\\
& \le 
\sup_{\substack{f \in \cK \\ \|e\|_p \le \zeta }}
\Big|\gamma(f) - \delta^{\rm opti} \Big(\La(f)+ \f{\eps}{\zeta}e \Big) \Big|
+ \sup_{\substack{f \in \cK \\ \|e\|_p \le \zeta }}
\Big|\delta^{\rm opti}(\La(f)+e) - \delta^{\rm opti}\Big( \La(f)+ \f{\eps}{\zeta}e \Big) \Big|\\
& \le e_\eps(\delta^{\rm opti})
+ \sup_{\|e\|_p \le \zeta} \sup_{\substack{a \in A \\ b \in B}} \Big| \sum_{k=1}^m c^{(a,b)}_k \Big( 1 - \f{\eps}{\zeta} \Big) e_k \Big) \Big|
= \inf_\delta e_\eps(\delta)
+ \sup_{\substack{a \in A \\ b \in B}} \| c^{(a,b)}\|_{p'} \Big( 1 - \f{\eps}{\zeta} \Big)  \zeta .
%\\
%& = \inf_\delta e_\eps(\delta) + \sup_{\substack{a \in A \\ b \in B}}  \| c^{(a,b)}\|_{p'} \;(\zeta - \eps).
\end{align*}
This indeed reduces to the desired inequality.
}

\section{Computational Realizations}
\label{SecComp}

The purpose of this section is to provide computational recipes for the practical construction of the optimal estimation functional appearing in Theorem~\ref{Thm4Max}
and, to a lesser extent, Theorem~\ref{ThmExt}.
At first, a rather abstract optimization program outputting the desired coefficient vectors will be presented.
Then, for the sake of implementation in a reproducible {\sc matlab} file\footnote{Accessible from the author's webpage or at \url{https://github.com/foucart/COR}.},
two situations of specific interest will be spelled out in separate subsections.

\subsection{Generic  convex optimization programs}

In most of this subsection and throughout the next two,
the nonlinear functional $\gamma: F \to \bR$ denotes a sup-linear functional,
i.e., it has the form $\gamma(f) = \sup_{i \in I} \gamma_i(f)$ with the $\gamma_i$'s being linear functionals.
It was established in \revX{Theorem~\ref{Thm4Max}} that the functionals $\delta: \bR^m \to \bR$ minimizing the estimation error $e(\delta) = \sup\{ |\gamma(f) - \delta(\La(f))|: \, f \in \cK \}$ can be sought among sup-affine functionals.
As a preliminary,
for a fixed such $\delta$,
one highlights how $e(\delta)$ can be computed by solving a convex optimization program.
The crucial point is that  $\nu := \delta \circ \La: F \to \bR$ is itself a sup-affine functional in this case.
The statement below features the standard simplex associated with the index set $I$,
defined when $I$ is finite by
$$
\cS^I := \bigg\{ \sigma \in \bR^I: \sigma_i \ge 0 \mbox{ for all } i \in I
\mbox{ and } \sum_{i \in I} \sigma_i = 1 \bigg\}.
$$

\bprop
\label{PropGWCEFixed}
Let $\cK$ be a convex set.
For a sup-linear functional $\gamma: F \to \bR$
and a sup-affine functional $\nu: F \to \bR$ given, for $f \in F$, by
$$
\gamma(f) = \sup_{i \in I} \gamma_i(f)
\qquad \mbox{and} \qquad
\nu(f) = \sup_{i \in I} \big( \nu_i(f) + b_i \big),
$$
\revvme{with the $\gamma_i, \nu_i \in F^*$ being linear functionals and the $b_i \in \bR$ being scalars},
one has 
\revvY{
\be
\label{CompGWCEFixed}
\sup_{f \in \cK} |\gamma(f) - \nu(f)|
 = 
\inf_{\substack{e \in \bR \\ \sigma^{(i_\star)}, \tau^{(i_\star)} \in \cS^I }} \; e 
\quad \mbox{s.to } 
\begin{cases}
\tnbig{\gamma_{i_\star} - \sum_{i \in I} \sigma^{(i_\star)}_i \nu_i }_\cK \le e + \sum_{i \in I} \sigma^{(i_\star)}_i b_i,\\
%& i_{\star} \in I,\\
\tnbig{\nu_{i_\star} - \sum_{i \in I} \tau^{(i_\star)}_i \gamma_i }_\cK \le e - b_{i_\star},
%& i_{\star} \in I.
\end{cases}
\ee
where the constraints over the scalar $e$ and the simplex-vectors $\sigma^{(i_\star)},\sigma^{(i_\star)}$
are understood to be imposed for all $i_\star \in I$.}
\eprop

\bpf
The supremum appearing on the left-hand side of \eqref{CompGWCEFixed}
is nothing else than the infimum value of $e$ subject to $|\gamma(f) - \nu(f)| \le e$ for all $f \in \cK$.
The latter decouples as the two constraints 
\be
\label{decouples}
\gamma(f) - \nu(f) \le e
\quad \mbox{for all } f \in \cK
\qquad \mbox{and} \qquad
\nu(f) - \gamma(f) \le e
\quad \mbox{for all } f \in \cK.
\ee
In view of $\gamma(f) = \sup_{i \in I} \gamma_i(f)$
and of 
$$
- \nu(f) 
= - \sup_{i \in I} \big( \nu_i(f) + b_i \big)
= \inf_{i \in I} \big( -\nu_i(f) - b_i \big)
= \inf_{\sigma \in \cS^I} \sum_{i \in I} \sigma_i \big( -\nu_i(f) - b_i \big),
$$
the first constraint in \eqref{decouples} is equivalent to a set of constraints indexed by $i_\star \in I$, namely
$$
\gamma_{i_\star}(f) + \inf_{\sigma \in \cS^I} \sum_{i \in I} \sigma_i \big( -\nu_i(f) - b_i \big) 
\le e 
\qquad \mbox{for all } f \in \cK. 
$$
Each of these individual constraints reads
$$
\sup_{f \in \cK} \inf_{\sigma \in \cS^I} \Big( \gamma_{i_\star}(f) + \sum_{i \in I} \sigma_i \big( -\nu_i(f) - b_i \big) \Big)
\le e.
$$
Since the inner expression depends affinely of $f \in \cK$,
which is a convex set,
and also affinely on $\sigma \in \cS^I$,
which is a convex and compact set,
von Neumann minimax theorem legitimizes the exchange of sup and inf,
so that the above constraint reads
$$
 \inf_{\sigma \in \cS^I} \sup_{f \in \cK} \Big( \gamma_{i_\star}(f) + \sum_{i \in I} \sigma_i \big( -\nu_i(f) - b_i \big) \Big)
\le e.
$$ 
All in all, the first constraint in \eqref{decouples} is equivalent to the existence of $\sigma^{(i_\star)} \in \cS^I$, $i_\star \in I$, such that
$$
\sup_{f \in \cK} \Big( \gamma_{i_\star}(f) - \sum_{i \in I} \sigma^{(i_\star)}_i \nu_i(f)  \Big)
\le e + \sum_{i \in I} \sigma^{(i_\star)}_i b_i,
$$ 
which is tautologically equivalent to the first set of convex constraints indexed by $i_\star$ appearing in~\eqref{CompGWCEFixed}. 
Likewise,
with very similar details left to the reader,
the second constraint in \eqref{decouples} is equivalent to the existence of $\tau^{(i_\star)} \in \cS^I$, $i_\star \in I$, obeying
the second set of convex constraints indexed by $i_\star$ appearing in~\eqref{CompGWCEFixed}. 
Incorporating all the $\sigma^{(i_\star)}$ and $\tau^{(i_\star)}$
as optimization variables,
one arrives at the announced convex optimization program.
\epf

Turning to the minimization of the estimation error $e(\delta)$ over all possible functionals $\delta: \bR^m \to \bR$,
or in fact only over all functionals of the form 
$\delta_c(y) = \sup_{i \in I} \big( c^{(i)}_0 + \sum_{k=1}^m c^{(i)}_k y_k \big)$, $y \in \bR^m$,
one notices that this error becomes
$$
e(\delta_c) = \sup_{f \in \cK} |\gamma(f) - \nu_c(f)|,
\qquad  \mbox{with} \quad
\nu_c := \delta_c \circ \Lambda : F \to \bR
\; \mbox{ being a sup-affine functional}.
$$ 
At a fixed $c \in (\bR \times \bR^m)^I$,
Proposition~\ref{PropGWCEFixed} just indicated how to compute $e(\delta_c)$ by solving a convex optimization program.
Regrettably,
a straightforward minimization over $c$ as well seems out of reach,
due to the presence of products $\sigma^{(i_\star)}_i c^{(i)}_k$ of optimization variables generated through the first set of constraints. 
An~alternative route,
pursued below,
consists in translating the constructive argument of Theorem~\ref{Thm4Max} into a manageable optimization program.

\bprop
\label{ThmConstruction}
Given a convex set $\cK$ containing the origin,
an optimal functional $\delta^{\rm opti}: \bR^m \to \bR$ for the estimation of a sup-linear functional $\gamma(f) = \sup_{i \in I} \gamma_i(f)$, $f \in F$,
is obtained as the sup-affine functional $\delta^{\rm opti}(y) = \sup_{i \in I} \big( \wh{c}^{(i)}_0  + \sum_{k=1}^m \wh{c}^{(i)}_k y_k \big)$, $y \in \bR^m$.
\revvme{The} `hat'-notation denotes solutions to the convex optimization program
\revvY{
\be
\label{OptProgGWCE}
\underset{\substack{e \in \bR, e',e'' \in \bR^I\\
c^{(i_\star)} \in \bR^m\\ 
\sigma^{(i_\star)} \in \cS^I}}
{\rm minimize}  \quad e
\qquad 
\mbox{s.to }
\begin{cases}
e'_{i_\star} + e''_{i_\star} \le 2e,\\
% & i_\star \in I,\\
\tnbig{ \gamma_{i_\star} - \sum_{k = 1}^m c^{(i_\star)}_k \la_k }_\cK \le e'_{i_\star}, \\
%& i_\star \in I,\\
\tnbig{ - \sum_{i \in I} \sigma^{(i_\star)}_i \gamma_i + \sum_{k = 1}^m c^{(i_\star)}_k \la_k }_\cK
\le e''_{i_\star}, 
\end{cases}
\ee
where the constraints on the scalars $e,e'_{i_\star},e''_{i_\star}$,
the vectors $c^{(i_\star)}$,
and the simplex-vectors $\sigma^{(i_\star)}$
are understood to be imposed for all $i_\star \in I$. 
}
Subsequently,
the remaining coefficients $\wh{c}^{(i)}_0 \in \bR$ are obtained as $\wh{c}^{(i_\star)}_0 = \wh{e} - \wh{e}''_{i_\star}$ for all $i_\star \in I$.
\eprop

\bpf
Keeping the proof of Theorem~\ref{Thm4Max} in mind,
recall that the key resided in \eqref{KEY} being satisfied for all $i_\star \in I$,
i.e., 
$$
\sup_{f' \in \cK}
\bigg( \gamma_{i_\star}(f') 
- \sum_{k = 1}^m c^{(i_\star)}_k \la_k(f') \bigg)
+
\sup_{f'' \in \cK}
\bigg(
-\gamma(f'') + \sum_{k = 1}^m c^{(i_\star)}_k \la_k(f'')
\bigg)
\le 2 e
$$
for some $c^{(i_\star) } \in \bR^m$ and for $e \in \bR$ as small as possible,
i.e., $e = e_\flat$.
Thus, introducing slack variables $e'_{i_\star}, e''_{i_\star} \in \bR$ bounding the above suprema,
the task consists in minimizing $e$ subject to $e'_{i_\star} + e''_{i_\star} \le 2 e$,
to $\sup_{f' \in \cK}
\big( \gamma_{i_\star}(f') 
- \sum_{k = 1}^m c^{(i_\star)}_k \la_k(f') \big) \le e'_{i_\star}$,
and to $\sup_{f'' \in \cK}
\big(
-\gamma(f'') + \sum_{k = 1}^m c^{(i_\star)}_k \la_k(f'')
\big) \le e''_{i_\star}$ for all $i_\star \in I$.
The first of these constraints are exactly the first constraints in \eqref{OptProgGWCE},
while the second ones reduce to the second constraints  in \eqref{OptProgGWCE} simply by the definition of the support function.
Considering now the third and last constraints, in view of 
$$
-\gamma(f'') = - \sup_{i \in I} \gamma_i(f'')
= \inf_{i \in I} (-\gamma_i(f'')) 
= \inf_{\sigma \in \cS^I}
 \sum_{i\in I} \sigma_i (-\gamma_i(f'')),
$$
one obtains
\begin{align*}
\sup_{f'' \in \cK}
\bigg(
-\gamma(f'') + \sum_{k = 1}^m c^{(i_\star)}_k \la_k(f'')
\bigg)
& = \sup_{f'' \in \cK}
\inf_{\sigma \in \cS^I}
\bigg( - \sum_{i \in I} \sigma_i \gamma_i(f'') + \sum_{k = 1}^m c^{(i_\star)}_k \la_k(f'')
\bigg) \\
& = \inf_{\sigma \in \cS^I} \;
\sup_{f'' \in \cK}
\bigg( - \sum_{i \in I} \sigma_i \gamma_i(f'') + \sum_{k = 1}^m c^{(i_\star)}_k \la_k(f'')
\bigg), 
\end{align*}
where the exchange of sup and inf relied on
von Neumann minimax theorem,
in the same way as for the proof of Proposition \ref{PropGWCEFixed}. 
As a result, the constraints can be rephrased as the existence of $\sigma^{(i_\star)} \in \cS^I$ such that
$$
\tnBig{ - \sum_{i \in I} \sigma^{(i_\star)}_i \gamma_i + \sum_{k = 1}^m c^{(i_\star)}_k \la_k }_\cK
\le e''_{i_\star}.
$$
Incorporating all the $\sigma^{(i_\star)} \in \cS^I$ as optimization variables,
one arrives at the third constraints in~\eqref{OptProgGWCE}.
All in all,
it has been justified that the convex optimization program \eqref{OptProgGWCE} outputs $e_\flat$ 
(i.e., the minimal error $e(\delta)$ over all possible estimation functionals $\delta$) as $\wh{e}$
and the coefficient vectors making the linear part of $\delta^{\rm opti}$ as the $\wh{c}^{(i_\star)} \in \bR^m$.
Finally, 
the announced expression for the remaining coefficients $\wh{c}^{(i_\star)}_0 \in \bR$ follows
from their very definition \eqref{DefC0}, written as
$$
\wh{c}^{(i_\star)}_0 = 
- \sup_{f'' \in \cK}
\bigg(
-\gamma(f'') + \sum_{k = 1}^m \wh{c}^{(i_\star)}_k \la_k(f'')
\bigg) + e_\flat,  
$$
and from the identification of the above supremum as $\wh{e}''_{i_\star}$ and of $e_\flat$ as $\wh{e}$. 
\epf

When estimating a mixed sup-inf-linear functional,
it is as also possible to translate the proof of Theorem~\ref{ThmExt} into a practical recipe 
for the construction of the optimal sup-inf-affine  functional.
The key resides again in fulfilling the inequalities \eqref{KEYInfSup},
which are treated as above.
More precisely, one proceeds in exactly the same way (save for the notation)
to recast the bound on the second summand in the left-hand side of of \eqref{KEYInfSup} by a slack variable $e''_{a_\star,b_\star}$,
and the argument will also apply when recasting a bound on the fist summand by a slack variable $e'_{a_\star,b_\star}$.
With details left out,
one arrives at the convex optimization program displayed below.

\bprop
\label{ThmConstruction2}
Under the setting of Theorem~\ref{ThmExt},
the optimal estimation functional $\delta^{\rm opti}: \bR^m \to \bR$ 
has the sup-inf-affine form $\delta^{\rm opti}(y) = \sup_{a \in A} \inf_{b \in B} \big( \wh{c}^{(a,b)}_0  + \sum_{k=1}^m \wh{c}^{(a,b)}_k y_k \big)$, $y \in \bR^m$.
\revvme{The} `hat'-notation denotes solutions to the convex optimization program
\revvY{
$$
%\label{OptProgGWCE}
\underset{\substack{e \in \bR, e',e'' \in \bR^{A \times B}\\
c^{(a_\star,b_\star)} \in \bR^m\\ 
\sigma^{(a_\star,b_\star)} \in \cS^{I_{a_\star}}\\
\tau^{(a_\star,b_\star)} \in \cS^{I_{b_\star}} }}
{\rm minimize}  \quad e
\qquad 
\mbox{s.to }
\begin{cases}
e'_{a_\star,b_\star} + e''_{a_\star,b_\star} \le 2e,\\ 
%& a_\star \in A, b_\star \in B,\\
\tnbig{ \sum_{i \in I_{a_\star}} \sigma^{(a_\star,b_\star)}_i \gamma_i - \sum_{k = 1}^m c^{(a_\star,b_\star)}_k \la_k }_\cK
\le e'_{a_\star,b_\star},\\
% & a_\star \in A, b_\star \in B,\\
\tnbig{ - \sum_{j \in J_{b_\star}} \tau^{(a_\star,b_\star)}_j \gamma_j + \sum_{k = 1}^m c^{(a_\star,b_\star)}_k \la_k }_\cK
\le e''_{a_\star,b_\star},
% & a_\star \in A, b_\star \in B.
\end{cases}
$$
where the constraints on the scalars $e,e'_{a_\star,b_\star},e''_{a_\star,b_\star}$,
the vectors $c^{(a_\star,b_\star)}$,
and the simplex-vectors $\sigma^{(a_\star,b_\star)},\tau^{(a_\star,b_\star)}$ are understood to be imposed for all $a_\star \in A, b_\star \in B$. 
}
Subsequently,
the coefficients $\wh{c}^{(a,b)}_0 \in \bR$ are obtained as $\wh{c}^{(a_\star,b_\star)}_0 = \wh{e} - \wh{e}''_{a_\star,b_\star}$ for all $a_\star \in A,  b_\star \in B$.
\eprop

Arguably, Propositions~\ref{ThmConstruction}
and~\ref{ThmConstruction2}
are still abstract statements.
Whether they really translate into practical constructions depends on the specific model set $\cK$ and its amenability to numerical computations.
Concerning Proposition~\ref{ThmConstruction2},
there is the added issue of the number of constraints,
which is $|A| \times |B|$.
In the example of the $\ell$th largest value among $d$ linear functionals,
one has
$$
|A| = |\{ I \in [1 \hspace{-.5mm} : \hspace{-.5mm}d]: |I| = \ell \}| = \binom{d}{\ell}
\quad \mbox{and} \quad
|B| = |\{ J \in [1 \hspace{-.5mm} : \hspace{-.5mm}d]: |J| = d + 1 - \ell \}| = \binom{d}{\ell-1},
$$
so $|A| \times |B|$ behaves exponentially in $d$ unless $\ell$ or $d-\ell$ is a small constant.
For this reason, one concentrates only on the estimation of the maximum of linear functionals when 
expounding on the two relevant examples presented in the next subsections.

\subsection{Polytopal models in $\ell_\infty$-spaces}

Working in $F = \ell_\infty^N$,
the convex model set containing the origin (but not necessarily symmetric) considered here is the polytope
$$
\cK^{\rm poly} := \big\{ f \in \ell_\infty^N: \langle a_\ell, f \rangle \le 1 \mbox{ for all } \ell \in L \big\}
$$ 
\revZ{described through prescribed vectors $a_\ell \in \ell_1^N$
indexed by a finite set $L$}.
Using duality in linear programming,
it is routine to verify that the support function of $\cK^{\rm poly}$, evaluated at a linear functional $\eta$ defined for $f \in \ell_\infty^N$ by $\eta(f) = \langle g, f \rangle$,
can be expressed as
\begin{align*}
\tn{\eta}_{\cK^{\rm poly}}
& = \sup_{f \in \bR^N} \langle g, f \rangle \quad
\mbox{s.to } \langle a_\ell, f \rangle \le 1
\mbox{ for all } \ell \in L\\
& = \inf_{s \in \bR^L} \sum_{\ell \in L} s_\ell
\quad
\mbox{s.to } s \ge 0 \mbox{ and } \sum_{\ell \in L} s_\ell a_\ell = g.
\end{align*}
Therefore,
an inequality constraint $\tn{\eta}_{\cK^{\rm poly}} \le \kappa$ found as a constraint in an optimization program can be rephrased as the existence of some $s \in \bR_{+}^L$ such that $\sum_{\ell \in L} s_\ell \le \kappa$ and $\sum_{\ell \in L} s_\ell a_\ell = g$,
and subsequently $s$ can be incorporated in the optimization variables.
For instance,
with $u_k \in \ell_1^N$ representing the linear functional $\la_k$ and
$w_i \in \ell_1^N$ representing the linear functional $\gamma_i$, so that
$$
\la_k(f) = \langle u_k, f \rangle
\quad \mbox{and} \quad
\gamma_i(f) = \langle w_i, f \rangle 
\qquad \mbox{for all } f \in \ell_\infty^N,
$$
the convex optimization program from Proposition~\ref{ThmConstruction} becomes:

\ul{\sl Linear program for the optimal estimation of $\gamma = \sup_{i \in I} \gamma_i$ with model $\cK^{\rm poly}$.}
\revvme{
$$
\underset{\substack{e \in \bR, e',e'' \in \bR^I\\
c^{(i_\star)} \in \bR^m\\ 
\sigma^{(i_\star)} \in \cS^I\\ s^{(i_\star)},t^{(i_\star)} \in \bR_+^L}}
{\rm minimize}  \quad e
\qquad 
\mbox{s.to }
%\begin{cases}
\begin{Bmatrix}
e'_{i_\star} + e''_{i_\star} \le 2e,
\; \sum_{\ell \in L} s^{(i_\star)}_\ell \le e'_{i_\star},
\; \sum_{\ell \in L} t^{(i_\star)}_\ell \le e''_{i_\star}\\
%, & i_\star \in I,\\
\sum_{\ell \in L} s^{(i_\star)}_\ell a_\ell =  w_{i_\star} - \sum_{k = 1}^m c^{(i_\star)}_k u_k\\
%, & i_\star \in I,\\
\sum_{\ell \in L} t^{(i_\star)}_\ell a_\ell =
- \sum_{i} \sigma^{(i_\star)}_i w_i + \sum_{k = 1}^m c^{(i_\star)}_k u_k
% , & i_\star \in I.
%\end{cases}
\end{Bmatrix}
\; \forall i_\star \in I.
$$

}

After solving this linear program,
one can construct an optimal estimation functional $\delta^{\rm opti}$ which has the generic sup-affine form
\be
\label{GenSupAff}
\delta_{b,z}(y) = \sup_{i \in I} \big( b_i + \langle z_i, y \rangle \big), 
\qquad y \in \bR^m.
\ee
For good measure, 
one can verify 
\revvme{(as done in the reproducible)}
that the minimal value of $e$ found above agrees with $e(\delta^{\rm opti})$
when the latter is computed by solving the convex program from Proposition~\ref{PropGWCEFixed}. 
In the present case, it becomes:

\ul{\sl Linear program to compute $e(\delta_{b,z})$ for a fixed sup-affine functional with model $\cK^{\rm poly}$.}
\revvme{
\be
\label{LP4Fixed}
\underset{\substack{e \in \bR \\ \sigma^{(i_\star)}, \tau^{(i_\star)} \in \cS^I\\ s^{(i_\star)},t^{(i_\star)} \in \bR_+^L }}{\rm minimize} \quad e 
\qquad \mbox{s.to } 
%\begin{cases}
\begin{Bmatrix}
\sum_{\ell \in L} s^{(i_\star)}_\ell \le e + \sum_{i \in I} \sigma^{(i_\star)}_i b_i,
\; \sum_{\ell \in L} t^{(i_\star)}_\ell \le e - b_{i_\star}\\
%, & i_{\star} \in I,\\
\sum_{\ell \in L} s^{(i_\star)}_\ell a_\ell = 
w_{i_\star} - \sum_{i \in I} \sigma^{(i_\star)}_i \La^* z_i \\
%, & i_{\star} \in I,\\
\sum_{\ell \in L} t^{(i_\star)}_\ell a_\ell =
\La^* z_{i_\star} - \sum_{i \in I} \tau^{(i_\star)}_i w_i\\
%, & i_{\star} \in I.
%\end{cases}
\end{Bmatrix}
\; \forall i_\star \in I.
\ee
}

%\vspace{5mm}

\brk
In $F = \ell_\infty^N$ and with a convex model,
it is somewhat known that there exists an affine map $\Delta^{\rm opti}: \bR^m \to \ell_\infty^N$
which is genuinely optimal for the estimation of \revvme{$\Id_F$},
i.e., for the full approximation problem.
This is essentially Sukharev's result applied componentwise.
Leaving out the details, this optimal map has the form
$$
\Delta^{\rm opti}:
y \in \bR^m
\mapsto
\bigg[
\wh{c}^{(0)} + \sum_{k=1}^m y_k \wh{c}^{(k)}
\bigg] \in \ell_\infty^N,
$$
where the vectors $\wh{c}^{(0)}, \wh{c}^{(1)},\ldots, \wh{c}^{(m)} \in \ell_\infty^N$ are solutions to the following linear program
(involving the standard basis $(v_1,\ldots,v_N)$ for $\bR^N$):
\revvme{
$$
\underset{\substack{e \in \bR \\ 
c^{(0)},c^{(k)} \in \bR^N\\ s^{(j)},t^{(j)} \in \bR_+^L }}{\rm minimize} \quad e 
\qquad \mbox{s.to } 
%\begin{cases}
\begin{Bmatrix}
\sum_{\ell \in L} s^{(j)}_\ell \le e + c^{(0)}_j,
\; \sum_{\ell \in L} t^{(j)}_\ell \le e - c^{(0)}_j\\
%, & j \in [1:N],\\
\sum_{\ell \in L} s^{(j)}_\ell a_\ell = 
v_j - \sum_{k=1}^m c^{(k)}_j u_k\\ 
%, & j \in [1:N],\\
\sum_{\ell \in L} t^{(j)}_\ell a_\ell =
-v_j + \sum_{k=1}^m c^{(k)}_j u_k
%,& j \in [1:N].
%end{cases}
\end{Bmatrix}
\; \forall j \in [1 \hspace{-.5mm} : \hspace{-.5mm} N]
$$ 
}As an initially guess, a natural candidate for the estimation of $\gamma = \sup_{i \in I} \gamma_i$ could have been the `plug-in' functional $\delta^{\rm plug} = \gamma \circ \Delta^{\rm opti}: \bR^m \to \bR$.
Since it has the generic sup-affine form \eqref{GenSupAff} with $b_i = \gamma_i(\wh{c}^{(0)})$ and $(z_{i})_{k} = \gamma_i(\wh{c}^{(k)})$,
the estimation error $e(\delta^{\rm plug})$ of this `plug-in' functional can be computed by solving the linear program \eqref{LP4Fixed}.
Numerical experiments available in the reproducible file confirm that this guess is not optimal, i.e., that $e(\delta^{\rm plug})$ is larger than $e(\delta^{\rm opti})$.
\erk

\subsection{Approximability models in Hilbert spaces}

Working now in a Hilbert space $F=H$,
the convex model set containing the origin (but not necessarily symmetric)
considered here is the set of elements that are approximated with prescribed accuracy $\eps > 0$ by an affine subspace $g+V$,
where $V$ is an $n$-dimensional linear subspace of $H$
and  $g \in H \setm V$ is a fixed element satisfying ${\rm dist}(g,V) < \eps$.
More precisely, one considers 
$$
\cK^{\rm appr} := \{ f \in H: {\rm dist}(f,g+V) \le \eps \}
= \{ f \in H: \|P_{V^\perp} (f-g)\| \le \eps \},
$$
where $P_{V^\perp}$ denotes the orthogonal projector onto the orthogonal complement of $V$.
Motivated by Optimal Recovery questions,
such approximability sets (with $g=0$) 
were introduced in \cite{BCDDPW}.
Based on an observation made (with $g=0$) in \cite[Subsection 3.3]{FouVV},
their support function evaluated at a linear functional $\eta$ can be expressed as
$$
\tn{\eta}_{\cK^{\rm appr}} = \eta(g) +
\begin{cases}
\eps \|\eta\|_{*}  & \mbox{if } \eta_{| V} = 0,\\
+\infty & \mbox{otherwise}. 
\end{cases}
$$
Consequently,
still using the notation
$u_k \in H$ for the Riesz representer of the linear functional $\la_k$
and $w_i \in H$ for the Riesz representers of the linear functional $\gamma_i$---in a reproducing kernel Hilbert space with kernel ${\sf k}$, 
note that $w_i = {\sf k}(\cdot, x^{(i)})$
if $\gamma_i$ is the evaluation at a point $x^{(i)}$---the convex optimization program from Proposition \ref{ThmConstruction} becomes:

\ul{\sl Second-order cone program for the optimal estimation of $\gamma = \sup_{i \in I} \gamma_i$ with model $\cK^{\rm appr}$.}
\revvme{
\begin{align*}
& \underset{\substack{e \in \bR, e',e'' \in \bR^I\\
c^{(i_\star)} \in \bR^m\\ 
\sigma^{(i_\star)} \in \cS^I}}
{\rm minimize}  \qquad e
\\
& 
%\quad \; \; \; 
\; \, \mbox{s.to }
%\begin{cases}
\begin{Bmatrix}
e'_{i_\star} + e''_{i_\star} \le 2e\\
%, & i_\star \in I,\\
w_{i_\star} - \sum_{k = 1}^m c^{(i_\star)}_k u_k \in V^\perp,
\; 
- \sum_{i \in I} \tau^{(i_\star)}_i w_i + \sum_{k = 1}^m c^{(i_\star)}_k u_k \in V^\perp\\
%, & i_\star \in I,\\
\eps \big\| w_{i_\star} - \sum_{k = 1}^m c^{(i_\star)}_k u_k \big\| \le e'_{i_\star} -  \langle w_{i_\star} - \sum_{k = 1}^m c^{(i_\star)}_k u_k, g \rangle\\
% &,  i_\star \in I,\\
\eps \big\| - \sum_{i \in I} \sigma^{(i_\star)}_i w_i + \sum_{k = 1}^m c^{(i_\star)}_k u_k \big\|
\le e''_{i_\star}
- \langle - \sum_{i \in I} \sigma^{(i_\star)}_i w_i + \sum_{k = 1}^m c^{(i_\star)}_k u_k , g \rangle
%, & i_\star \in I.
%\end{cases}
\end{Bmatrix}
\; \forall  i_\star \in I.
\end{align*}
}

As in the previous subsection,
one can construct an optimal estimation functional $\delta^{\rm opti}$ of the generic sup-affine form \eqref{GenSupAff} after solving the above second-order cone program.
And again, for good measure, 
one can verify
\revvme{(as done in the reproducible)}
 that the minimal value of $e$ just found  agrees with $e(\delta^{\rm opti})$
when the latter is computed by solving the convex program from Proposition~\ref{PropGWCEFixed}. 
In the present case, it becomes:

\ul{\sl Second-order cone program to compute $e(\delta_{b,z})$ for a fixed sup-affine functional  with model $\cK^{\rm appr}$.}
\revvme{
\begin{align}
\label{SOCP4Fixed}
& \underset{\substack{e \in \bR \\ \sigma^{(i_\star)}, \tau^{(i_\star)} \in \cS^I }}{\rm minimize} \quad e
\\
\nonumber 
& 
%\quad \; \; \; \; \, 
\; \; \mbox{s.to } 
%\begin{cases}
\begin{Bmatrix}
w_{i_\star} - \sum_{i \in I} \sigma^{(i_\star)}_i \La^* z_i \in V^\perp,
\; 
\La^* z_{i_\star} - \sum_{i \in I} \tau^{(i_\star)}_i w_i \in V^\perp\\
%, & i_\star \in I,\\
\eps \| w_{i_\star} - \sum_{i \in I} \sigma^{(i_\star)}_i \La^* z_i\| 
\le  e + \sum_{i \in I} \sigma^{(i_\star)}_i b_i
- \langle w_{i_\star} - \sum_{i \in I} \sigma^{(i_\star)}_i \La^* z_i , g \rangle\\
%, & i_{\star} \in I,\\
\eps \| \La^* z_{i_\star} - \sum_{i \in I} \tau^{(i_\star)}_i w_i \|
\le  e - b_{i_\star} - \langle \La^* z_{i_\star} - \sum_{i \in I} \tau^{(i_\star)}_i w_i , g \rangle
%, & i_{\star} \in I.
%\end{cases}
\end{Bmatrix}
\; \forall  i_\star \in I.
\end{align}
}

%\vspace{5mm}

\brk
In $F=H$ and with the model set $\cK^{\rm appr}$,
there also exists an affine map $\Delta^{\rm opti}: \bR^m \to H$ which is genuinely optimal for the estimation of \revvme{$ \Id_F$}.
Specifically,
it was established in \cite{BCDDPW} that,
$$
\Delta^{\rm cheb}: y \in \bR^m \mapsto 
\bigg[ \underset{f \in H}{\argmin} \|P_{V^\perp} f \|
\quad \mbox{s.to }
\La(f) = y \bigg]
\in H
$$
is an optimal\footnote{In fact, it was even shown that $\Delta^{\rm cheb}(y)$ is the Chebyshev center for the set $\cK^{\rm appr} \cap \La^{-1}(\{y\})$---i.e., the center of the smallest ball containing the set---for any $y \in \La(\cK^{\rm appr})$, which is stronger than the notion of optimality examined so far.} recovery map in the case $g = 0$.
It then easily follows that 
$$
\Delta^{\rm opti}:
y \in \bR^m \mapsto g + \Delta^{\rm cheb}(y - \La(g)) \in H
$$ 
is an optimal recovery map with $g \not= 0$.
It is also affine,
since $\Delta^{\rm cheb}$ is a linear map which admits the expression (not immediate from above, see \cite[Section 10.2]{MPDSE} for details)
$$
\Delta^{\rm cheb}(y)
= \sum_{k=1}^m a_k u_k + \sum_{\ell=1}^n b_\ell v_\ell, 
\qquad
a = \big[ I_m - C (C^\top C)^{-1} C^\top \big] y,
\quad
b = \big[ (C^\top C)^{-1} C^\top \big] y,
$$
where the $u_k$ are assumed to be orthonormal,
the $v_\ell$ form a basis for $V$,
and $C \in \bR^{m \times n}$ is the cross-gramian matrix with entries $\langle u_k,v_\ell \rangle$.
Given the above,
an initial candidate for the estimation of $\gamma = \sup_{i \in I} \gamma_i$ is the `plug-in' functional $\delta^{\rm plug} = \gamma \circ \Delta^{\rm opt}: \bR^m \to \bR$.
The latter has the generic sup-affine form \eqref{GenSupAff} with $b_i = \langle w_i,g - \Delta^{\rm cheb}\La g \rangle$ and $z_i = (\Delta^{\rm cheb})^* w_i$,
so the estimation error $e(\delta^{\rm plug})$ of this `plug-in' functional can be computed by solving the second-order cone program~\eqref{SOCP4Fixed}.
Numerical~experiments available in the reproducible file hold a surprise:
on all the examples tested,
the `plug-in' functional was optimal,
i.e., $e(\delta^{\rm plug})$ was equal to $e(\delta^{\rm opti})$.
In Hilbert spaces, 
for model sets of the form $\{ f \in H: \|T f\| \le 1 \}$ with $T$ not necessarily equal to $(1/\eps)P_{V^\perp}$,
it was known (see \cite[Lemma 6]{FouHen}) that `plug-in' maps are indeed optimal for the estimation of linear maps,
but the indication that this phenomenon could remain true for the estimation of some nonlinear maps is quite striking.
\erk

\end{document}